\documentclass[a4paper,english]{amsart}
\usepackage{setspace}
\usepackage{amssymb}
\IfFileExists{url.sty}{\usepackage{url}}
                      {\newcommand{\url}{\texttt}}

\makeatletter
 \theoremstyle{plain}
\newtheorem{thm}{Theorem}[section]
  \theoremstyle{plain}
  \newtheorem{cor}[thm]{Corollary}
  \theoremstyle{definition}
  \newtheorem{defn}[thm]{Definition}
  \theoremstyle{plain}
  \newtheorem{lem}[thm]{Lemma}
  \theoremstyle{plain}
  \newtheorem{prop}[thm]{Proposition}

\usepackage{babel}
\makeatother
\begin{document}

\title{Ekedahl-Oort strata in the supersingular locus}

\author{Maarten Hoeve}

\begin{abstract}
We give a description of the individual Ekedahl-Oort strata contained
in the supersingular locus in terms of Deligne-Lusztig varieties,
refining a result of Harashita.
\end{abstract}

\address{Korteweg-de Vries Instituut, Universiteit van Amsterdam, Plantage
Muidergracht 24, 1018 TV Amsterdam, The Netherlands.}

\email{mhoeve@science.uva.nl}

\maketitle

\section{Introduction}

Fix a prime $p$ and consider the moduli stack $\mathcal{A}_{g}$ of
principally polarized abelian varieties of dimension $g$ over
$\mathbb{F}_{p}$. In \cite{Oo01} Ekedahl and Oort defined a
stratification of $\mathcal{A}_{g}$, the Ekedahl-Oort
stratification. Harashita showed that certain unions of strata in
the supersingular locus are isomorphic to Deligne-Lusztig
varieties (see \cite{Ha07}). Here we show that individual strata
are isomorphic to a finer kind of Deligne-Lusztig varieties.

The difference between these results is clearest in the index
sets. Van der Geer and Moonen showed how to index Ekedahl-Oort
strata using the Weyl group\[ W_{g}=\{
w\in\mathrm{Aut}(\{1,\dots,2g\})\,|\, w(2g+1-i)=2g+1-w(i)\}\] of
the symplectic group $\mathrm{Sp}_{2g}$: there is an open stratum
$\mathcal{S}_{w}$ for each $w$ in \[ ^{I}W_{g}=\{ w\in W_{g}\,|\,
w^{-1}(1)<\dots<w^{-1}(g)\},\] see \cite{Mo01}. Equivalently, we
can use cosets: if $W_{g,I}$ is the subgroup generated by the
permutations $(i,i+1)(2g-i,2g+1-i)$ for $i=1,\dots,g-1$, then the
natural map $^{I}W_{g}\to W_{g,I}\backslash W_{g}$ is a bijection.

We now summarize Harashita's description of the Ekedahl-Oort
strata in the supersingular locus. The stratum attached to
$w\in{}^{I}W_{g}$ is in this locus if and only if $w$ is in the
subgroup \[ W_{g}^{[c]}=\{ w\in{}W_{g}\,|\, w(i)=i\textrm{ for
}i=1,\dots,g-c\}\] for some $c\leq g/2$. This subgroup is
isomorphic to $W_{c}$ via
\[
\mathfrak{r}:{}W_{g}^{[c]}\stackrel{\sim}{\to}{}W_{c},\quad\mathfrak{r}(w)(i)=w(g-c+i)-(g-c).\]
Write $^{I}W_{g}^{[c]}=W_{g}^{[c]}\cap{}^{I}W_{g}$ and
$^{I}W_{g}^{(c)}={}^{I}W_{g}^{[c]}-{}^{I}W_{g}^{[c-1]}$. Note that
$\mathfrak{r}$ maps $^{I}W_{g}^{[c]}$ to $^{I}W_{c}$.

Fix a $c\leq g/2$ and a supersingular elliptic curve $E$ over
$\mathbb{F}_{p^{2}}$ and let $\Lambda_{g,c}$ be the set of
isomorphism classes of polarizations on $E^{g}$ with kernel
isomorphic to $\alpha_{p}^{2c}$. For any $\mu\in\Lambda_{g,c}$ the
isogenies from $(E^{g},\mu)$ to principally polarized abelian
varieties correspond one-to-one with maximal isotropic subspaces
in a $2c$-dimensional symplectic vector space (see section
\ref{sec:X-as-moduli-space}). Let $X_{0}$ be the variety over
$\mathbb{F}_{p^{2}}$ that parameterizes both objects.

On the one hand we get a morphism $i_{\mu}\colon
X_{0}\to\mathcal{A}_{g,\mathbb{F}_{p^{2}}}$ that sends an isogeny
to its target and an action of $\mathrm{Aut}(E^{g},\mu)$ on
$X_{0}$ by precomposition. The morphism $i_{\mu}$ factors through
the quotient stack $[X_{0}/\mathrm{Aut}(E^{g},\mu)]$.

On the other hand we get Deligne-Lusztig varieties $X_{0}[\alpha]$
in $X_{0}$. For $\alpha\in W_{c,I}\backslash W_{c}/W_{c,I}$ the
variety $X_{0}[\alpha]$ consists of all parabolic subspaces $U$
such that $U$ and $U^{(p^{2})}$ (or rather their stabilizers) are
in relative position $\alpha$.

\begin{thm}
\label{thm:Harashita-Main-Theorem}(\cite{Ha07}, main theorem) For
each double coset $\alpha=W_{c,I}w'W_{c,I}$ with
$w'\in{}^{I}W_{g}^{[c]}$ the morphism of stacks over
$\mathbb{F}_{p}$\begin{eqnarray*}
\coprod_{\mu\in\Lambda_{g,c}}[X_{0}[\alpha]/\mathrm{Aut}(E^{g},\mu)]
& \to & \bigcup_{w\in{}^{I}W_{g}\textrm{ s.t.
}\mathfrak{r}(w)\in\alpha}\mathcal{S}_{w}\end{eqnarray*} induced
by the $i_{\mu}$ is finite and surjective. If $w'$ is in
$^{I}W_{g}^{(c)}$, then this morphism is a bijection on geometric
points.
\end{thm}
Unfortunately, this theorem only describes unions of strata. In
this paper we refine it to individual strata by looking at a finer
kind of Deligne-Lusztig varieties. Inspired by ideas of Moonen and
Wedhorn in \cite{MW04}, we consider not just the relative position
of a subspace $U$ with $U^{(p^{2})}$, but also that of their
refinements. We get varieties $X_{0}(w')$ indexed by $^{I}W_{c}$
(instead of $W_{c,I}\backslash W_{c}/W_{c,I}$). With these we can
refine theorem \ref{thm:Harashita-Main-Theorem} to the following.

\begin{thm}
\label{thm:Refinement-Harashita}For each $w$ in $^{I}W_{g}^{(c)}$
the morphism of stacks over $\mathbb{F}_{p}$\[
\coprod_{\mu\in\Lambda_{g,c}}[X_{0}(\mathfrak{r}(w))/\mathrm{Aut}(E^{g},\mu)]\to\mathcal{S}_{w},\]
induced by the $i_{\mu}$ is an isomorphism.
\end{thm}
An Ekedahl-Oort stratum is reducible if it is contained in the
supersingular locus, at least if $p$ is large enough (see
\cite{Ha07} corollary 3.5.3) and irreducible otherwise (see
\cite{EG06} theorem 11.5). In section
\ref{sec:Number-of-components} we show that
$X_{0}(\mathfrak{r}(w))$ is irreducible for
$w\in{}^{I}W_{g}^{(c)}$. So theorem \ref{thm:Refinement-Harashita}
gives the exact number of components of the supersingular strata.

\begin{cor}
\label{cor:Number-Of-Components}For $w\in{}^{I}W_{g}^{(c)}$ the
Ekedahl-Oort stratum $\mathcal{S}_{w}$ has $\#\Lambda_{g,c}$
irreducible components.
\end{cor}
Here it is crucial that $w$ is in $^{I}W_{g}^{(c)}$ and not just
in $^{I}W_{g}^{[c]}$. The number $\#\Lambda_{g,c}$ is a class
number, see \cite{Ha07} section 3.2 and 3.5.

This description of the EO-strata in the supersingular locus bears
a striking resemblance to the description of Kottwitz-Rapoport
strata in the supersingular locus in a recent paper \cite{GY08} by G\"ortz
and Yu. It would be interesting to further investigate the
relation between Ekedahl-Oort and Kottwitz-Rapoport strata.

\subsubsection{Conventions}

If $S$ is a scheme over $\mathbb{F}_{p}$ and $\mathcal{F}$ is an
$\mathcal{O}_{S}$-module, then we denote with
$\mathcal{F}^{(p^{r})}=\mathcal{O}_{S}\otimes_{\mathcal{O}_{S},\mathrm{Frob}^{r}}\mathcal{F}$
the pull-back  by the $r$-th power of the absolute Frobenius.
Further, if $S\to T$ is a morphism of schemes and $X$ a
$T$-scheme, we write $X_{S}$ for $S\times_{T}X$.

\subsubsection{Acknowledgements}

I thank my advisors Ben Moonen and Gerard van der Geer for their
help. I thank Torsten Wedhorn for pointing out a mistake in an
earlier version.

\section{\label{sec:Deligne-Lusztig-varieties}Deligne-Lusztig varieties}

We can improve Harashita's theorem, because we use finer
Deligne-Lusztig varieties. Originally Deligne and Lusztig defined
their varieties for Borel subgroups (see \cite{DL76}). We can
generalize them to parabolic subgroups in two ways. If we use
exactly the same definition, we get what we will call coarse
Deligne-Lusztig varieties. By theorem
\ref{thm:Harashita-Main-Theorem} they correspond to unions of
strata. If we use a more subtle definition, pioneered by Lusztig
and B\'edard (see section 1.2 of \cite{Lu03} and \cite{Be85}),
then we get fine Deligne-Lusztig varieties. Theorem
\ref{thm:Refinement-Harashita} says that they correspond to
individual strata.

In this section we gather some definitions and results on
Deligne-Lusztig varieties. Although we will only need them for
$\mathrm{Sp}_{2c}$, we give them for more general groups.

\subsection{Coarse Deligne-Lusztig varieties}

Let $k$ be an algebraic closure of the finite field with $q$
elements $\mathbb{F}_{q}$. Suppose that $G$ is a reductive
connected algebraic group over $k$, obtained by extension of
scalars from $G_{0}$ over $\mathbb{F}_{q}$. Write $F\colon G\to G$
for the corresponding Frobenius morphism.

Let $W$ be the Weyl group of $G$. The Frobenius $F$ acts on $W$.
Let $S\subset W$ be the set of reflections in simple roots. We
denote with $\mathcal{P}_{I}$ the variety over $k$ of parabolic
subgroups in $G$ of type $I\subset S$. The group $G$ acts on
$\mathcal{P}_{I}$ by conjugation. So it also acts on
$\mathcal{P}_{I}\times\mathcal{P}_{J}$ for all $I,J\subset S$.

Let $W_{I}$ is the subgroup of $W$ generated by the reflections in
the roots in $I\subset S$. The $G$-orbits in
$\mathcal{P}_{I}\times\mathcal{P}_{J}$ are in bijection with
$W_{I}\backslash W/W_{J}$, see \cite{Be85} II lemma 7. Suppose $P$
and $Q$ are two parabolic subgroups, of types $I$ and $J$
respectively. Then we say that they are in \emph{relative
position} $w\in W_{I}\backslash W/W_{J}$ if the point
$(P,Q)\in\mathcal{P}_{I}\times\mathcal{P}_{J}$ is in the $G$-orbit
corresponding to $w$. We write $\mathrm{relpos}(P,Q)=w$.

Each double coset in $W_{I}\backslash W/W_{J}$ contains a unique
element of minimal length. Denote the set of such elements with
$^{I}W^{J}\subset W$, so that the quotient map $^{I}W^{J}\to
W_{I}\backslash W/W_{J}$ is a bijection. We will often see this
bijection as an identification. In particular, we will speak of
parabolic subgroups in relative position $w\in{}^{I}W^{J}$.

\begin{defn}
The \emph{coarse Deligne-Lusztig variety} $\mathcal{P}_{I}[w]$
attached to $w\in{}^{I}W^{F(I)}$ is the locally closed subscheme
of $\mathcal{P}_{I}$ consisting of all parabolic subgroups $P$
such that $P$ and $F(P)$ are in relative position $w$.
\end{defn}
The orbit in $\mathcal{P}_{I}\times\mathcal{P}_{F(I)}$
corresponding to $w\in{}^{I}W^{F(I)}$ is smooth of dimension
$l(w)+\dim(\mathcal{P}_{I\cap F(I)})$, where $l$ is the length
function. The variety $\mathcal{P}_{I}[w]$ is the intersection of
this orbit with the graph of the Frobenius morphism. Since the
intersection is transversal, we get the following (compare
\cite{DL76} section 1.3).

\begin{lem}
\label{lem:Basic-prop-coarse-var}The variety $\mathcal{P}_{I}[w]$
is smooth and purely of dimension $l(w)+\dim(\mathcal{P}_{I\cap
F(I)})-\dim(\mathcal{P}_{I})$.
\end{lem}
In particular, if $F(I)=I$, then $\mathcal{P}_{I}[w]$ has
dimension $l(w)$. There is also the following result of Bonnaf\'e
and Rouquier \cite{BR06} on irreducibility.

\begin{thm}
\label{lem:Irreducibility-coarse-var}A variety
$\mathcal{P}_{I}[w]$ is reducible if and only if $W_{I}w$ is
contained in a proper subgroup of the form $W_{J}$ for some
$J\subset S$ that is $F$-stable.
\end{thm}

Using this one can determine the number of irreducible components
of any $\mathcal{P}_{I}[w]$, see \cite{GY08} corollary 5.3.

\subsection{\label{sub:Fine-Deligne-Lusztig-varieties}Fine Deligne-Lusztig varieties}

Keep the notation from the previous section. To get finer
varieties, we look not just at the relative position of $P$ and
$F(P)$, but also of their refinements. Given two parabolic
subgroups $P$ and $Q$ of $G$, the \emph{refinement of $P$ with
respect to $Q$} is\[ \mathrm{Ref}_{Q}(P)=(P\cap
Q)U_{P}=U_{P}(P\cap Q),\] where $U_{P}$ is the unipotent radical
of $P$. This is again a parabolic subgroup and it is contained in
$P$. If $P$ is of type $I$ and $Q$ of type $J$ and they are in
relative position $w$, then $\mathrm{Ref}_{Q}(P)$ is of type
$I\cap{}^{w}J$.

Suppose $I$ is a subset of the set of simple roots. Given a
sequence $\mathbf{u}=(u_{0},u_{1},\dots)$ of elements of $W$,
define a sequence of subsets $I_{n}\subset I$ by $I_{0}=I$ and
$I_{n+1}=I_{n}\cap{}^{u_{n}}F(I_{n})$. Let $\mathcal{T}(I)$ be the
set of sequences $\mathbf{u}$, such that\[
u_{n}\in{}^{I_{n}}W^{F(I_{n})}\quad\textrm{and}\quad u_{n+1}\in
W_{I_{n+1}}u_{n}W_{F(I_{n})}.\] Then we have the following
description of $\mathcal{T}(I)$, see \cite{Be85} I proposition 9.

\begin{prop}
\label{pro:Biject-sequences-JW}Each sequence $\mathbf{u}$ in
$\mathcal{T}(I)$ stabilizes to some $u_{\infty}$, i.e. satisfies
$u_{n}=u_{n+1}=\dots=u_{\infty}$ for $n$ large enough. The map
$\mathcal{T}(I)\to W$ which sends $\mathbf{u}$ to $u_{\infty}$
induces a bijection from $\mathcal{T}(I)$ to $^{I}W$.
\end{prop}
\begin{defn}
\label{def:Fine-Deligne-Lustig-var}The \emph{fine Deligne-Lustig
variety} $\mathcal{P}_{I}(\mathbf{u})$ attached to a sequence
$\mathbf{u}$ in $\mathcal{T}(I)$ is the locally closed subscheme
of $\mathcal{P}_{I}$ consisting of all parabolic subgroups $P$
such that if we define\[ P_{0}=P\quad\textrm{and}\quad
P_{n+1}=\mathrm{Ref}_{F(P_{n})}(P_{n})\] then $P_{n}$ and
$F(P_{n})$ are in relative position $u_{n}$.
\end{defn}

See \cite{Lu03} 1.3 and 1.4 for some examples of fine Deligne-Lusztig varieties.

It follows from proposition \ref{pro:Biject-sequences-JW} that we
can also index the fine Deligne-Lusztig varieties by $^{I}W$. So
we will often speak of $\mathcal{P}_{I}(w)$ for $w\in{}^{I}W$.

Most questions about the fine varieties can be reduced to ones
about coarse varieties, using the following proposition of
B\'edard, see \cite{Be85} II proposition 12.

\begin{prop}
The morphism
$\mathcal{P}_{I_{0}}(u_{0},u_{1},\dots)\to\mathcal{P}_{I_{1}}(u_{1},u_{2},\dots)$
which sends $P$ to $\mathrm{Ref}_{F(P)}(P)$ is an isomorphism.
\end{prop}
Iterating the above proposition until we get to $u_{\infty}$ we
get the following.

\begin{cor}
There is an isomorphism
$\mathcal{P}_{I}(\mathbf{u})\cong\mathcal{P}_{I_{\infty}}[u_{\infty}]$,
where $I_{\infty}=\bigcap I_{n}$.
\end{cor}
Now we combine this with lemmas \ref{lem:Basic-prop-coarse-var}
and \ref{lem:Irreducibility-coarse-var}.

\begin{cor}
\label{cor:properties-fine-DL-var} The variety
$\mathcal{P}_{I}(\mathbf{u})$ is smooth and purely of dimension
$l(u_{\infty})+\dim(\mathcal{P}_{I_{\infty}\cap
F(I_{\infty})})-\dim(\mathcal{P}_{I_{\infty}})$. It is reducible
if and only if $W_{I_{\infty}}u_{\infty}$ is contained in a proper
$F$-stable standard parabolic subgroup of $W$.
\end{cor}
When $G=\mathrm{Sp}_{2c}$, the Frobenius $F$ is the identity on
$W$, so $\mathcal{P}_{I}(\mathbf{u})$ has dimension
$l(u_{\infty})$.

The action of $G$ on $\mathcal{P}_{I}$ by conjugation restricts to
an action of $G_{0}(\mathbb{F}_{q})$. For $g\in
G_{0}(\mathbb{F}_{q})$ we have $F(g)=g$. Hence,
$F(gPg^{-1})=gF(P)g^{-1}$. The relative position of two parabolic
subgroups is unchanged if you conjugate them by the same element.
So both the coarse and fine Deligne-Lusztig varieties are stable
under the $G_{0}(\mathbb{F}_{q})$-action.

\subsubsection{An equivalent definition}

It will be convenient to have the following description of the
fine Deligne-Lusztig varieties.

\begin{prop}
The variety $\mathcal{P}_{I}(\mathbf{u})$ consists of all
parabolics $P$ in $\mathcal{P}_{I}$ such that if we set\[
P_{0}'=P\quad\textrm{and}\quad
P_{n+1}'=\mathrm{Ref}_{F(P_{n}')}(P),\] then $P$ and
$F(P_{\infty}')$ are in relative position $u_{\infty}$.
\end{prop}
This description differs in two ways from the definition of
$\mathcal{P}_{I}(\mathbf{u})$. First of all it uses $P_{n}'$
instead of $P_{n}$. That $P_{n}'=P_{n}$ follows by induction from
the following lemma with $P=P$, $Q'=F(P_{n}')$ and
$Q=F(P_{n-1}')$.

\begin{lem}
Suppose $Q'\subset Q$ and $P$ are parabolic subgroups. Then \[
\mathrm{Ref}_{Q'}(\mathrm{Ref}_{Q}(P))=\mathrm{Ref}_{Q'}(P).\]

\end{lem}
\begin{proof}
The left hand side is $U_{P}((U_{P}(P\cap Q)\cap Q')$. So it is
sufficient to show that $(U_{P}(P\cap Q))\cap Q'=P\cap Q'$. On the
one hand, $P\cap Q\subset U_{P}(P\cap Q)$, gives $P\cap Q'=P\cap
Q\cap Q'\subset(U_{P}(P\cap Q))\cap Q'$. On the other hand,
$P\supseteq U_{P}(P\cap Q)$ gives $P\cap Q'\supseteq(U_{P}(P\cap
Q))\cap Q'$.
\end{proof}
Second of all the condition that $P$ and $F(P_{\infty}')$ are in
relative position $u_{\infty}$ is equivalent to the condition that
$P_{\infty}$ and $F(P_{\infty})$ are in relative position $w$ by
the following lemma with $Q=F(P_{\infty})$.

\begin{lem}
For any two parabolic subgroups $P$ and $Q$ one has\[
\mathrm{relpos}(P,Q)=\mathrm{relpos}(\mathrm{Ref}_{Q}(P),Q)\in
W.\]

\end{lem}
\begin{proof}
This is the case $Z=P$ in lemma 3.2(c) in \cite{Lu03}.
\end{proof}
Since the sequence $\mathbf{u}$ is determined by $u_{\infty}$
(lemma \ref{pro:Biject-sequences-JW}), the fact that $P_{\infty}$
and $F(P_{\infty})$ are in relative position $u_{\infty}$ is
equivalent with $P$ being in $\mathcal{P}_{I}(\mathbf{u})$. This
proves the proposition.

\subsection{\label{sub:DL-var-for-Sp}The case that $G_{0}=\mathrm{Sp}_{2c}$}

Let $L_{0}$ be a $2c$-dimensional vector space over
$\mathbb{F}_{q}$ with a symplectic form and let
$G_{0}=\mathrm{Sp}(L_{0})$ be the symplectic group of $L_{0}$. A
partial flag in $L:=k\otimes L_{0}$ is a collection of subspaces
$\mathcal{C}$ that is totally ordered by the inclusion. There is a
bijection between such flags and parabolic subgroups of $G$ by
sending a flag $\mathcal{C}$ to its stabilizer
$\mathrm{stab}(\mathcal{C})$.

The Weyl group of $G$ is the group $W_{c}$ from the introduction.
For $w\in W_{c}$ we define\[ r_{w}(i,j)=\#\{
a\in\{1,\dots,i\}\,|\, w(a)\leq j\}.\] Using elementary linear
algebra one can prove the following.

\begin{lem}
\label{lem:Relpos-in-flags}Suppose $\mathcal{C}$ and $\mathcal{D}$
are two flags in $L_{0}$ and let $P$ and $Q$ be their respective
stabilizers. Let $I$ be the type of $P$ and $J$ that of $Q$. Then
$P$ and $Q$ are in relative position $w\in{}^{I}W^{J}$ if and only
if\[ \dim(C\cap D)=r_{w}(\dim(C),\dim(D))\quad\textrm{for all
}C\in\mathcal{C}\textrm{ and }D\in\mathcal{D}.\]

\end{lem}
If $P$ is the stabilizer of $\mathcal{C}=\{0=C_{0}\subset
C_{1}\subset\dots\subset C_{r}=L\}$ and $Q$ that of
$\mathcal{D}=\{0=D_{0}\subset D_{1}\subset\dots\subset D_{s}=L\}$,
then $\mathrm{Ref}_{Q}(P)$ is the stabilizer of the flag
consisting of\[ (C_{i+1}\cap D_{j})+C_{i}\textrm{ for
}i=1,\dots,r\textrm{ and }j=1,\dots,s.\] Denote this flag with
$\mathrm{Ref}_{\mathcal{D}}(\mathcal{C})$. Also, note that if $P$
is the stabilizer of $\mathcal{C}$, then $F(P)$ is the stabilizer
of $\mathcal{C}^{(q)}=\{ C_{i}^{(q)}\}$.

\section{\label{sec:X-as-moduli-space}A moduli space of isogenies}

In this section we construct an isomorphism between the variety of
isogenies from $(E^{g},\mu)$ to principally polarized abelian
varieties and a variety of maximal isotropic subspaces in a
symplectic vector space. This is a direct generalization of the
construction of the families of abelian surfaces considered by
Moret-Bailly in \cite{MB81}.

\subsection{\label{sub:Dieudonn=E9-modules}Dieudonn\'e modules}

Suppose that $S$ is a scheme over $\mathbb{F}_{p}$. Let
$\hat{\mathrm{CW}}$ be the formal group of Witt-covectors over
$\mathbb{Z}$, see \cite{Fo77}, chapter II 1.5. For a formal
$p$-group scheme $G$ over $S$ let
\[
M(G)=\mathcal{H}om(G,\hat{\mathrm{CW}}_{S}),\] be the sheaf on $S$
whose sections over $U\subset S$ are the homomorphisms
$G|_{U}\to\hat{\mathrm{CW}}_{U}$. There is a natural action of the
Witt-vectors $W$ on $\hat{\mathrm{CW}}$. This makes $M(G)$ into a
$W(\mathcal{O}_{S})$-module. The Frobenius and Verschiebung on $G$
give homomorphisms $F\colon M(G)^{(p)}\to M(G)$ and $V\colon
M(G)\to M(G)^{(p)}$ respectively (where
$M(G)^{(p)}=W(\mathcal{O}_{S})\otimes_{W(\mathcal{O}_{S}),\mathrm{Frob}}M(G)$).
By abuse of notation we will sometimes write $M(G)$ for the triple
$(M(G),F,V)$.

Over the spectrum of a perfect field $M(G)$ is the classical
contravariant Dieudonn\'e module of $G$ as defined by Fontaine in
\cite{Fo77}. The functor $G\rightsquigarrow M(G)$ is exact and
gives an equivalence from the category of formal $p$-groups to the
category of modules over a certain ring, see loc. cit. theorems 1
and 2 in the introduction. Over other schemes $S$, this is no
longer true in general, but we will still call $M(G)$ the
\emph{Dieudonn\'e module} of $G$.

If $G$ is annihilated by $p$, multiplication by $p$ is zero on
$M(G)$. Hence, $M(G)$ is actually on $\mathcal{O}_{S}$-module. If
$G$ is annihilated by $V$, then
$M(G)=\mathcal{H}om(G,\hat{\mathbb{G}_{a}})$ by proposition III,
3.2 in \cite{Fo77}. It follows from \cite{Jo93} proposition 2.2
that for any $S$ the functor $M$ is an equivalence from the
categories of group schemes that are annihilated by $V$ to that of
$\mathcal{O}_{S}$-modules $M$ with a homomorphism $F:M^{(p)}\to
M$.

If $\pi:A\to S$ is an abelian scheme over $S$, then there is an
isomorphism\[ M(A[p])^{(p)}\cong
H_{\mathrm{dR}}^{1}(A/S)\stackrel{\mathrm{def}}{=}\mathbb{R}^{1}\pi_{*}(\Omega_{A/S}^{*})\]
by \cite{BBM82} theorem 4.2.14.

\subsection{The moduli space}

Like in \cite{LO98} section 1.2, we fix a supersingular elliptic
curve $E'$ over $\mathbb{F}_{p}$ such that the Frobenius satisfies
$F^{2}+p=0$. All the endomorphisms of $E'$ are defined over
$\mathbb{F}_{p^{2}}$, so it is convenient to work with
$E=E'_{\mathbb{F}_{p^{2}}}$.

Fix a $g\geq1$. Every polarization of $E^{g}$ is defined over
$\mathbb{F}_{p^{2}}$, since $E^{g}$ is isomorphic to its own dual,
so a polarizations can be seen as an endomorphism. For $c\leq
g/2$, let $\Lambda_{g,c}$ be the set of isomorphism classes of
polarizations $\mu$ on $E^{g}$ whose kernel is isomorphic to
$\alpha_{p}^{2c}$.

Now fix a polarization $\mu\in\Lambda_{g,c}$. If we put
$\overline{M}_{0}=M(E^{g}[p])$, then $\mu$ gives a morphism of
Dieudonn\'e modules
$\mu^{*}\colon\overline{M}_{0}\to\overline{M}_{0}^{t}$ , where
$\overline{M}_{0}^{t}$ is the dual of $\overline{M}_{0}$. Let
$L_{0}$ be the cokernel of $\mu^{*}$. Then $L_{0}$ is the
Dieudonn\'e module of the kernel of $\mu$ and so it is isomorphic
to $\mathbb{F}_{p^{2}}^{2c}$ with $F$ and $V$ both the zero map.
The polarization $\mu$ induces a pairing $e^{\mu}$ on its kernel
and this gives a symplectic form on $L_{0}$.

Let $S$ be any scheme over $\mathbb{F}_{p^{2}}$. We want to study
the set $I_{\mu}(S)$ of isomorphism classes of isogenies \[
\rho\colon(E_{S}^{g},\mu_{S})\to(A,\lambda),\] where $A$ is an
abelian scheme over $S$ and $\lambda$ is a principal polarization
on $A$. Note that $\rho$ gives a morphism\[ \rho^{*}\colon
M(A[p])\to
M(E_{S}^{g}[p])=\mathcal{O}_{S}\otimes_{\mathbb{F}_{p^{2}}}\overline{M}_{0}\]
of Dieudonn\'e modules on $S$.

\begin{lem}
\label{lem:Descr-moduli-space}The map that sends an isogeny $\rho$
to the image of the composition of $\rho^{*}$ with
$\mathcal{O}_{S}\otimes_{\mathbb{F}_{p^{2}}}\overline{M}_{0}\to\mathcal{O}_{S}\otimes_{\mathbb{F}_{p^{2}}}L_{0}$,
gives a bijection from $I_{\mu}(S)$ to the set of isotropic subbundles
of rank $c$ in $\mathcal{O}_{S}\otimes_{\mathbb{F}_{p^{2}}}L_{0}$.
\end{lem}

\begin{proof}
An isogeny $\rho$ is determined up to isomorphism by its kernel.
Because $\mu_{S}=\rho^{*}\lambda$, we have
$\ker(\rho)\subset\ker(\mu_{S})$. Since
$\ker(\mu_{S})\cong\alpha_{p}^{2c}$ is annihilated by $V$, the
inclusion $\ker(\rho)\subset\ker(\mu_{S})$ is determined by the
induced morphism\[
\mathcal{O}_{S}\otimes_{\mathbb{F}_{p^{2}}}L_{0}=M(\ker(\mu_{S}))\to
M(\ker(\rho^{*}))=\mathrm{coker}(\rho^{*})\] on Dieudonn\'e
modules. This morphism is again determined by its kernel, which is
exactly the image of $\rho^{*}$ composed with $\overline{M}_{0}\to
L_{0}$.

Let us determine which subbundles of
$\mathcal{O}_{S}\otimes_{\mathbb{F}_{p^{2}}}L_{0}$ one can get in
this way. Any subbundle of
$\mathcal{O}_{S}\otimes_{\mathbb{F}_{p^{2}}}L_{0}$ gives by
Dieudonn\'e theory a subgroup scheme $G\subset\ker(\mu_{S})$. The
polarization $\mu$ descends to $E_{S}^{g}/G$ if and only if $G$ is
isotropic with respect to $e^{\mu}$ and if it descends, say to
$\lambda$, then the kernel of $\lambda$ is $G^{\bot}/G$.
Therefore, one gets all maximal isotropic subspaces of
$\mathcal{O}_{S}\otimes_{\mathbb{F}_{p^{2}}}L_{0}$.
\end{proof}
Let $X_{0}$ be the variety over $\mathbb{F}_{p^{2}}$ that
parameterizes isotropic subspaces of rank $c$ in $L_{0}$. Note that as
a space it depends only on $c$. The
lemma shows that $X_{0}$ is a fine moduli space for the functor
$S\rightsquigarrow I(S)$. In particular, there is a universal
isogeny\[
\rho_{X_{0}}\colon(E_{X_{0}}^{g},\mu_{X_{0}})\to(A_{X_{0}},\lambda_{X_{0}})\]
over $X_{0}$. The principally polarized abelian scheme
$(A_{X_{0}},\lambda_{X_{0}})$ induces a morphism \[ i_{\mu}\colon
X_{0}\to\mathcal{A}_{g,\mathbb{F}_{p^{2}}}\] over
$\mathbb{F}_{p^{2}}$.

Let $k$ be an algebraic closure of $\mathbb{F}_{p}$. Write $X$ for
the base change of $X_{0}$ to $k$. Using the bijection between
flags and parabolics from section \ref{sub:DL-var-for-Sp}, we can
see $X$ as a variety of parabolics $\mathcal{P}_{I}$ for some $I$
(the type of the stabilizer of a maximal isotropic subgroup).
In particular there are Deligne-Lusztig varieties $X(w)$ in $X$
for $w\in{}^{I}W_{c}$.

\subsection{\label{sub:The-action-of-Aut}The action of $\mathrm{Aut}(E^{g},\mu)$}

The automorphism group $\Gamma_{\mu}:=\mathrm{Aut}(E^{g},\mu)$
acts on $I(S)$ by precomposition: $\gamma\in\Gamma_{\mu}$ sends an
isogeny $\rho\colon(E_{S}^{g},\mu_{S})\to(A,\lambda)$ to
$\rho\circ\gamma_{S}$. The action is functorial in $S$, so it
induces an action of $\Gamma_{\mu}$ on $X_{0}$.

Equivalently, we can define this action as follows. The action of
$\Gamma_{\mu}$ on $E^{g}$ induces actions on $\overline{M}_{0}$
and $\overline{M}_{0}^{t}$ such that $\mu^{*}$ is equivariant.
This gives an action of $\Gamma_{\mu}$ on the cokernel $L_{0}$
which respects the symplectic pairing, i.e. a homomorphism
$\Gamma_{\mu}\to\mathrm{Sp}(L_{0})(\mathbb{F}_{p^{2}})$. Via the
action of this last group on $X_{0}$, we get an action of
$\Gamma_{\mu}$ on $X_{0}$. Since the Deligne-Lusztig varieties are
$\mathrm{Sp}(L_{0})(\mathbb{F}_{p^{2}})$-stable, they are also
$\Gamma_{\mu}$-stable.

\section{\label{sec:The-EO-stratification}The EO-stratification on $X$}

In this section we pull the EO-stratification back to $X$ by the
morphism $i_{\mu}\colon X\to(\mathcal{A}_{g})_{k}$. Keep the
notation from the previous section. Suppose $x$ is a $k$-valued
point of $X$, corresponding to an isogeny\[
\rho\colon(E_{k}^{g},\lambda_{k})\to(A,\lambda).\] To see in which
EO-stratum $i_{\mu}(x)$ is, we must express the Dieudonn\'e module
$\overline{N}$ of $A[p]$ in terms of the subspace $U\subset
L:=k\otimes L_{0}$ attached to $x$. This is done using the
filtration constructed below.

\subsection{A filtration of the Dieudonn\'e modules}

It will be convenient to work with the $p^{\infty}$ Dieudonn\'e
modules. Let $M_{0}$ be the Dieudonn\'e module of
$E^{g}[p^{\infty}]$ and let $N$ be that of $A[p^{\infty}]$. These
are modules over the Witt vectors $W(\mathbb{F}_{p^{2}})$ and
$W(k)$ respectively. Put
$M=W(k)\otimes_{W(\mathbb{F}_{p^{2}})}M_{0}$. Then $\mu$ and
$\rho$ induce \emph{}injective morphisms\[ \mu^{*}\colon
M_{0}^{t}\hookrightarrow M_{0}\quad\textrm{and}\quad\rho^{*}\colon
N\hookrightarrow M.\] Since over a perfect field the Dieudonn\'e
functor is exact, we have $\overline{M}_{0}=M_{0}/pM_{0}$ and
$\overline{N}=N/pN$.

Define a descending filtration $\mathrm{Fil}^{*}$ on $M_{0}$ by
\[
\mathrm{Fil}^{i}(M_{0}) =
\begin{cases}
M_{0} & \textrm{for }i\leq0,\\
\mathrm{im}(\mu^{*}) & \textrm{for }i=1,\\
F((\mathrm{Fil}^{i-2})^{(p)}) & \textrm{for }i=2,3,\\
p\,\mathrm{Fil}^{i-4} & \textrm{for }i\geq4.
\end{cases}
\]
Then
$\mathrm{gr}^{0}(M_{0})=\overline{M}_{0}/\mathrm{im}(\mu^{*})=L_{0}$,
since the image of $\mu^{*}$ contains $pM_{0}$. For brevity we
write $K_{0}=\mathrm{gr}^{1}(M_{0})$ and $K=k\otimes K_{0}$ (this
space plays no significant role). Because of the inherent
periodicity of the filtration, the graded modules
are\begin{equation} \mathrm{gr}^{i+4j}(M_{0})\cong\begin{cases}
L_{0} & \textrm{for }i=0,\\
K_{0} & \textrm{for }i=1,\\
L_{0}^{(p)} & \textrm{for }i=2,\\
K_{0}^{(p)} & \textrm{for
}i=3,\end{cases}\label{eq:Graded-Modules-E}\end{equation} for each
$j\geq0$. The isomorphisms are induced by multiplication by
$-p^{j}$ for $i=0,1$ and by applying $-p^{j}F$ for $i=2,3$.

Since $F^{2}+p=0$ for our choice of $E$, we have $F=-V$ on $M_{0}$
and $V^{2}=-FV=-p$. So by construction
$V(\mathrm{Fil}^{i}M_{0})=\mathrm{Fil}^{i+2}(M_{0})^{(p)}$. The
isomorphisms in (\ref{eq:Graded-Modules-E}) are chosen so that the
graded maps
$\mathrm{gr}^{i}(M_{0})\to\mathrm{gr}^{i+2}(M_{0}^{(p)})$ induced
by $V$ are the identity. Note that $L_{0}^{(p^{2})}=L_{0}$, since
it is a vector space over $\mathbb{F}_{p^{2}}$.

The filtration $\mathrm{Fil}^{*}(M_{0})$ on $M_{0}$ is the main
tool in our analysis of the EO-strata on $X$ . By pulling it back
or taking it modulo $p$ we get filtrations on other spaces. We
denote them with $\mathrm{Fil}^{*}(\dots)$ with the space in
brackets, for instance $\mathrm{Fil}^{*}(M_{0}^{(p)})$ and
$\mathrm{Fil}^{*}(\overline{M}_{0})$. We now focus on the
filtration $\mathrm{Fil}^{*}(\overline{N})$, obtained by pulling
$\mathrm{Fil}^{*}(M_{0})$ back by $\rho^{*}$ and taking it modulo
$p$. It allows us to express $\overline{N}$ in terms of $L$ and
$U$.

\begin{lem}
\label{lem:Graded-Modules-A}The graded modules of $\overline{N}$
are\[ \mathrm{gr}^{i}(\overline{N})=\begin{cases}
U & \textrm{for }i=0,\\
k\otimes\mathrm{gr}^{i}(M_{0}) & \textrm{for }i=1,2,3,\\
L/U & \textrm{for }i=4,\\
0 & \textrm{otherwise,}\end{cases}\] where $U\subset L$ is the
subspace corresponding to $x$.
\end{lem}
\begin{proof}
We see $\mu^{*}$ and $\rho^{*}$ as inclusions: $M^{t}\subset
N\subset M$. Then $\mathrm{gr}^{0}(\overline{N})=N/M^{t}$ and this
quotient is $U$ by the construction of the universal isogeny on
$X$ in lemma \ref{lem:Descr-moduli-space}. Also, \[
\mathrm{gr}^{4}(\overline{N})=pM/pN\cong M/N=L/U.\] Because
$pN\subset pM=\mathrm{Fil}^{4}$ the graded modules for $i=1,2,3$
are unaltered (except for the extension of scalars to $k$).
\end{proof}
Again we have
$V(\mathrm{Fil}^{i}(\overline{N}))\subset\mathrm{Fil}^{i+2}(\overline{N})^{(p)}$.
So $V$ induces maps\[
\mathrm{gr}^{i}(V):\mathrm{gr}^{i}(\overline{N})\to\mathrm{gr}^{i+2}(\overline{N})^{(p)}.\]
Since on $M_{0}$ the graded maps induced by $V$ are the identity,
these maps are as follows.

\begin{lem}
\label{lem:Graded-Maps}Under the isomorphisms in lemma
\ref{lem:Graded-Modules-A}:
\begin{enumerate}
\item $\mathrm{gr}^{0}(V)$ is the inclusion $U\subset L$;
\item $\mathrm{gr}^{1}(V)$ is the identity on $K$;
\item $\mathrm{gr}^{2}(V)$ is the canonical surjection $L\to L/U$ and
\item $\mathrm{gr}^{i}(V)$ is $0$ for all other $i$.
\end{enumerate}
\end{lem}
A submodule $H\subset\mathrm{gr}^{i}(\overline{N})$ gives a
submodule of $\overline{N}$ by pulling it back via the projection
$\mathrm{pr}_{i}:\mathrm{Fil}^{i}\to\mathrm{gr}^{i}$. For such
submodules it is easy to calculate the pull-back by $F$ and $V$.

\begin{lem}
\label{lem:Graded-Pull-Back}For every $i$ and
$H\subset\mathrm{gr}^{i+2}(\overline{N})$ we have\begin{align*}
V^{-1}(\mathrm{pr}_{i+2}^{-1}(H)^{(p)}) & =\mathrm{pr}_{i}^{-1}(\mathrm{gr}^{i}(V)^{-1}(H^{(p)})),\\
F^{-1}(\mathrm{pr}_{i+2}^{-1}(H)) &
=\mathrm{pr}_{i}^{-1}(\mathrm{gr}^{i}(F)^{-1}(H)).\end{align*}

\end{lem}
\begin{proof}
This follows from the relation
$\mathrm{gr}^{i-2}(V)\circ\mathrm{pr}_{i-2}=\mathrm{pr}_{i}\circ
V$ defining $\mathrm{gr}^{i}(V)$ and the same relation for $F$.
\end{proof}
In particular this can be applied to the subspace $U=0$ of
$\mathrm{gr}^{4}(\overline{N})=L/U$ to get the following.

\begin{cor}
\label{cor:ker(F)-and-ker(V)}For $\overline{N}$, we have
$\ker(F)=\mathrm{pr}_{2}^{-1}(U)$ and
$\ker(V)=\mathrm{pr}_{2}^{-1}(U^{(p)})$.
\end{cor}

\subsection{Calculating the EO-type}

Now we prove the following proposition.

\begin{prop}
\label{pro:Pull-back-EO-to-X}The reduced underlying subscheme of
the pull-back\[
i_{\mu}^{-1}\mathcal{S}_{w}=X\times_{\mathcal{A}_{g,k}}\mathcal{S}_{w}\]
 of the (open) Ekedahl-Oort stratum $\mathcal{S}_{w}$ by $i_{\mu}\colon X\to\mathcal{A}_{g,k}$
is the Deligne-Lustig variety $X(\mathfrak{r}(w))$ if
$w\in{}^{I}W_{g}^{[c]}$. Otherwise it is empty.
\end{prop}
We must show that $i_{\mu}(x)$ is in $\mathcal{S}_{w}$ if and only
if $x$ is in $X(\mathfrak{r}(w))$. We do this by constructing
flags on $\overline{N}$ from flags on $L$.

\begin{defn}
Given symplectic flags $\mathcal{D}$ and $\mathcal{D}'$ in $L$
that contain $U$, define a flag
$\mathcal{E}(\mathcal{D},\mathcal{D}')$ in $\overline{N}$ as
follows. On $\mathrm{gr}^{0}(\overline{N})=U$ take all subspaces
in $\mathcal{D}$ contained in $U$, on $\mathrm{gr}^{4}=L/U$ take
those containing $U$ and on $\mathrm{gr}^{2}=L^{(p)}$ take the
flag $(\mathcal{D}')^{(p)}$. Finally pull everything back by
$\mathrm{pr}_{i}:\mathrm{Fil}^{i}\to\mathrm{gr}^{i}$.
\end{defn}
Remember the construction of the \emph{canonical flag} (see
\cite{Oo01} section 5, in particular lemma 5.2, 4): start with the
flag $\mathcal{C}_{0}=\{\ker(V)\}$ and create $\mathcal{C}_{i+1}$
from $\mathcal{C}_{i}$ by adding $V^{-1}(C)$ and
$V^{-1}(C)^{\bot}$ for $C\in\mathcal{C}_{i}$. The flag
$\mathcal{C}_{\infty}$ to which this sequence stabilizes is the
canonical flag.

Let $P$ be the stabilizer of $U$ and let $\mathcal{D}_{i}$ (resp.
$\mathcal{D}_{\infty}$) be the flag corresponding to the parabolic
subgroup $P_{i}$ (resp. $P_{\infty}$) in definition
\ref{def:Fine-Deligne-Lustig-var}. The flags $\mathcal{C}_{i}$ and
$\mathcal{D}_{i}$ are then related as follows.

\begin{prop}
For all $i$ \[
\mathcal{C}_{2i}=\mathcal{E}(\mathcal{D}_{i},\mathcal{D}_{i})\quad\textrm{and}\quad\mathcal{C}_{2i+1}=\mathcal{E}(\mathcal{D}_{i+1},\mathcal{D}_{i}).\]

\end{prop}
\begin{proof}
Use induction on $i$. The case $i=0$ follows from corollary
\ref{cor:ker(F)-and-ker(V)}. Suppose now that
$\mathcal{C}_{2i}=\mathcal{E}(\mathcal{D}_{i},\mathcal{D}_{i})$.
To calculate $\mathcal{C}_{2i+1}$, we need to pull back
$\mathcal{C}_{2i}^{(p)}$ by $V:\overline{N}\to\overline{N}^{(p)}$.
By lemma \ref{lem:Graded-Pull-Back}, we can do this on the graded
modules. Only in three cases is $\mathrm{gr}^{i}(V)$ non-zero:
\begin{enumerate}
\item $\mathrm{gr}^{0}(V)\colon U\to L=L^{(p^{2})}$ is the inclusion. The
restriction of $\mathcal{C}_{2i}^{(p)}$ to
$(\mathrm{gr}^{2})^{(p)}=L^{(p^{2})}$ is
$\mathcal{D}_{i}^{(p^{2})}$. So the pull-back by
$\mathrm{gr}^{0}(V)$ consists of $U\cap D^{(p^{2})}$ for
$D\in\mathcal{D}_{i}$, i.e. the subspaces of
$\mathrm{Ref}_{\mathcal{D}_{i}^{(p^{2})}}\mathcal{D}_{0}=\mathcal{D}_{i+1}$
contained in $U$.
\item $\mathrm{gr}^{1}(V)\colon K\to K$ is the identity. But the restriction
of $\mathcal{C}_{2i}^{(p)}$ to $(\mathrm{gr}^{3})^{(p)}=K$ is
empty. Hence, so it its pull-back.
\item $\mathrm{gr}^{2}(V)\colon L^{(p)}\to L^{(p)}/U^{(p)}$ is the quotient
map. The restriction of $\mathcal{C}_{2i}^{(p)}$ to
$(\mathrm{gr}^{4})^{(p)}=L^{(p)}/U^{(p)}$ consists of the
subspaces in $\mathcal{D}_{i}^{(p)}$ that contain $U^{(p)}$. The
pull-back to $L^{(p)}$ consists of the same subspaces.
\end{enumerate}
After adding the orthogonal complements, we get
$\mathcal{E}(\mathcal{D}_{i+1},\mathcal{D}_{i})$. The other case
is similar and we omit the details
\end{proof}
\begin{cor}
The canonical flag is
$\mathcal{C}_{\infty}=\mathcal{E}(\mathcal{D}_{\infty},\mathcal{D}_{\infty})$.
\end{cor}
Define $\psi_{w}(i)=i-r_{w}(g,i)$ with $r_{w}$ as in section
\ref{sub:DL-var-for-Sp}. Then $i_{\mu}(x)$ is in $\mathcal{S}_{w}$
if and only if $F(C^{(p)})$ has dimension $\psi_{w}(\dim(C))$ for
all $C\in\mathcal{C}_{\infty}$. This is equivalent to\[
\dim(\ker(F)\cap
C^{(p)})=\dim(C)-\psi_{w}(\dim(C))=r_{w}(g,\mathrm{dim}(C)).\]

By corollary \ref{cor:ker(F)-and-ker(V)} the kernel of $F$
contains $\mathrm{Fil}^{3}$ and is contained in
$\mathrm{Fil}^{2}$. So the intersection with a pull-back of
$D^{(p)}$ in $(\mathrm{gr}^{4})^{(p)}$, for
$D\in\mathcal{D}_{\infty}$, has dimension $\dim(D)$. The
intersection with a pull-back of the $D^{(p)}$in
$(\mathrm{gr}^{0})^{(p)}$ has dimension $g$. Hence, $i_{\mu}(x)$
is in $\mathcal{S}_{w}$ for some $w\in{}^{I}W_{g}^{[c]}$.

Now look at the pull-backs of the $D^{(p^{2})}$ in
$(\mathrm{gr}^{2})^{(p)}$. We have \[
\dim(\ker(F)\cap(p_{2}^{-1}D^{(p^{2})}))=g-c+\dim(U\cap
D^{(p^{2})}).\] One easily checks that for $w\in W_{g}^{[c]}$\[
r_{w}(g-c+i,g-c+j)=g-c+r_{\mathfrak{r}(w)}(i,j).\] So for
$w\in{}^{I}W_{g}^{[c]}$
\begin{align*}
i_{\mu}(x)\in\mathcal{S}_{w} &
\iff \dim(\ker(F)\cap C^{(p)})=r_{w}(g,\mathrm{dim}(C))\textrm{ for all }C\in\mathcal{C}_{\infty}\\
& \iff\dim(U\cap
D^{(p^{2})})=r_{\mathfrak{r}(w)}(c,\dim(D))\textrm{ for all } D
\in \mathcal{D}_{\infty}.
\end{align*}
By lemma \ref{lem:Relpos-in-flags} the last statement is
equivalent with $P_{0}$ and $P_{\infty}$ being in relative
position $\mathfrak{r}(w)$. This proves proposition
\ref{pro:Pull-back-EO-to-X}.

\section{The differential of $i_{\mu}$}

In this section we calculate the cotangent map of
$i_{\mu}:X\to(\mathcal{A}_{g})_{k}$ and show that it is
surjective.

We keep the notation of section \ref{sec:The-EO-stratification}.
In particular $X$ is the variety parameterizing $c$-dimensional
isotropic subspaces in a $2c$-dimensional vector space $L$. Write
$\mathcal{L}=\mathcal{O}_{X}\otimes_{k}L$ and let
$\mathcal{U}\subset\mathcal{L}$ be the universal subbundle on $X$.
Then we have the following well-known result.

\begin{lem}
\label{lem:Cotang-bundle-Grassm}The composition \[
\mathcal{U}\hookrightarrow\mathcal{L}\stackrel{d\otimes1}{\to}\Omega_{X}\otimes\mathcal{L}\to\Omega_{X}\otimes(\mathcal{L}/\mathcal{U})\]
is $\mathcal{O}_{X}$-linear. It induces an isomorphism
$\mathrm{Sym}^{2}(\mathcal{U})\to\Omega_{X}$.
\end{lem}
The universal isogeny $\rho_{X}:E_{X}^{g}\to A_{X}$ induces a
homomorphism\[ \rho_{X}^{*}:H_{\mathrm{dR}}^{1}(A_{X}/X)\to
H_{\mathrm{dR}}^{1}(E_{X}^{g}/X)\] which is horizontal with
respect to the Gauss-Manin connections on both sides. Note that \[
H_{\mathrm{dR}}^{1}(E_{X}^{g}/X)=\mathcal{O}_{X}\otimes_{\mathbb{F}_{p^{2}}}H_{\mathrm{dR}}^{1}(E^{g}/\mathbb{F}_{p^{2}})\cong\mathcal{O}_{X}\otimes_{\mathbb{F}_{p^{2}}}\overline{M}_{0}^{(p)}\]
(see section \ref{sub:Dieudonn=E9-modules}). With respect to these
isomorphisms, the Gauss-Manin connection is just $d\otimes1$.
Write $\mathcal{H}$ for $H_{\mathrm{dR}}^{1}(A_{X}/X)$ and let
$\mathrm{Fil}^{i}(\mathcal{H})$ be the pull-back of
$\mathcal{O}_{X}\otimes\mathrm{Fil}^{i}(\overline{M}_{0})^{(p)}$
by $\rho_{X}^{*}$.

\begin{lem}
\label{lem:GM-conn-and-filtrations}The subspaces
$\mathrm{Fil}^{i}(\mathcal{H})$ are horizontal with respect to the
Gauss-Manin connection $\nabla$, i.e.
$\nabla(\mathrm{Fil}^{i}(\mathcal{H}))\subset\Omega_{X}\otimes\mathrm{Fil}^{i}(\mathcal{H})$.
\end{lem}
\begin{proof}
Since Gauss-Manin connection on $H_{\mathrm{dR}}^{1}(E_{X}^{g}/X)$
is $d\otimes1$, the subspaces
$\mathcal{O}_{X}\otimes\mathrm{Fil}^{i}(\overline{M}_{0})^{(p)s}$are
horizontal. Since $\rho_{X}^{*}$ is horizontal, the same for their
pull-backs $\mathrm{Fil}^{i}(\mathcal{H})$.
\end{proof}
Let $\mathcal{E}=\pi_{*}\Omega_{A_{X}/X}\subset\mathcal{H}$ be the
Hodge bundle. Then $i_{\mu}^{*}\Omega_{(\mathcal{A}_{g})_{k}}$ is
isomorphic to $\mathrm{Sym}^{2}(\mathcal{E})$. Like in lemma
\ref{lem:Cotang-bundle-Grassm}, the composition \begin{equation}
\mathcal{E}\hookrightarrow\mathcal{H}\stackrel{\nabla}{\to}\Omega_{X}\otimes\mathcal{H}\to\Omega_{X}\otimes(\mathcal{H}/\mathcal{E})\label{eq:KS-map}\end{equation}
gives an $\mathcal{O}_{X}$-module homomorphism
$\mathrm{Sym}^{2}(\mathcal{E})\to\Omega_{X}$. This is exactly the
cotangent map of $i_{\mu}:X\to(\mathcal{A}_{g})_{k}$.

Note that $\rho_{X}^{*}$ induces an isomorphism\[
\mathrm{gr}^{2}(\mathcal{H})\stackrel{\rho^{*}}{\to}\mathcal{O}_{X}\otimes\mathrm{gr}^{2}(\overline{M}_{0}^{(p)})=\mathcal{O}_{X}\otimes
L=\mathcal{L}.\] In particular, we can see $\mathcal{U}$ as a
submodule of $\mathrm{gr}^{2}(\mathcal{H})$. For a closed point
$x$ of $X$ we have $k(x)\otimes\mathcal{E}=\ker(F)$ in
$\overline{N}^{(p)}=H_{\mathrm{dR}}^{1}(A/k)$ (notation as in
section \ref{sec:The-EO-stratification}). So it follows from
corollary \ref{cor:ker(F)-and-ker(V)} that $\mathcal{E}$ is the
pull-back of $\mathcal{U}$ by
$\mathrm{pr}_{2}:\mathrm{Fil}^{2}(\mathcal{H})\to\mathrm{gr}^{2}(\mathcal{H})$.

\begin{prop}
Under the isomorphisms
$\Omega_{X}\cong\mathrm{Sym}^{2}(\mathcal{U})$ (lemma
\ref{lem:Cotang-bundle-Grassm}) and
$i_{\mu}^{*}\Omega_{(\mathcal{A}_{g})_{k}}\cong\mathrm{Sym}^{2}(\mathcal{E})$
the cotangent map
$i_{\mu}^{*}\Omega_{\mathcal{A}_{g}}\to\Omega_{X}$ of $i_{\mu}$ is
just\[
\mathrm{Sym}^{2}(\mathcal{E})\to\mathrm{Sym}^{2}(\mathcal{E}/\mathrm{Fil}^{2}(\mathcal{H}))\stackrel{\rho_{X}^{*}}{\cong}\mathrm{Sym}^{2}(\mathcal{U}).\]

\end{prop}
\begin{proof}
By lemma \ref{lem:GM-conn-and-filtrations} the map
(\ref{eq:KS-map}) factors as\[
\mathcal{E}\to\mathcal{E}/\mathrm{Fil}^{3}(\mathcal{H})\stackrel{h}{\to}\Omega_{X}\otimes(\mathrm{Fil}^{2}(\mathcal{H})/\mathcal{E})\hookrightarrow\Omega_{X}\otimes(\mathcal{H}/\mathcal{E}).\]
Now $\rho_{X}^{*}$ map $\mathcal{E}/\mathrm{Fil}^{3}(\mathcal{H})$
isomorphically to $\mathcal{U}$ and
$\mathrm{Fil}^{2}(\mathcal{H})/\mathcal{E}$ isomorphically to
$\mathcal{L}/\mathcal{U}$. Since $\rho_{X}^{*}$ is horizontal and
on $\mathcal{O}_{X}\otimes\overline{M}_{0}^{(p)}$ the connection
is just $d\otimes1$, under $\rho_{X}^{*}$ the map $h$ becomes the
same map as the one in lemma \ref{lem:Cotang-bundle-Grassm}.
\end{proof}
In particular the cotangent map is surjective.

\section{Proof of the main theorem \ref{thm:Refinement-Harashita}}

To prove the main theorem, we need some facts about the
EO-stratification. Fix a $w$ in ${}^{I}W_{g}^{[c]}$. Oort showed
that the open Ekedahl-Oort stratum $\mathcal{S}_{w}$ is purely of
dimension $l(w)$ (see \cite{Oo01} theorem 1.2; the formulation in
terms of lengths is due to Moonen, \cite{Mo03}). Also it is smooth
by \cite{EG06} corollary 8.4 (this is implicit in \cite{Oo01}
where the tangent spaces are calculated).

Consider the composition \[
j_{\mu}:X(\mathfrak{r}(w))\hookrightarrow
i_{\mu}^{-1}\mathcal{S}_{w}\stackrel{i_{\mu}}{\to}\mathcal{S}_{w}.\]
It is proper, since the first map is a closed immersion and the
second one comes from $i_{\mu}\colon X\to(\mathcal{A}_{g})_{k}$.
Because both sides are smooth of dimension $l(w)$ and $i_{\mu}$ is
surjective on cotangent spaces, $j_{\mu}$ is \'etale. Since
$j_{\mu}$ is quasi-compact, it is quasi-finite and, hence, finite
\'etale.

The disjoint union of the $j_{\mu}$\[
\coprod_{\mu\in\Lambda_{g,c}}j_{\mu}\colon\coprod_{\mu\in\Lambda_{g,c}}X(\mathfrak{r}(w))\to\mathcal{S}_{w}.\]
is finite \'etale and also surjective by theorem
\ref{thm:Harashita-Main-Theorem}. Each $j_{\mu}$ is defined over
$\mathbb{F}_{p^{2}}$. However,
$\sigma\in\mathrm{Gal}(\mathbb{F}_{p^{2}}/\mathbb{F}_{p})$ sends
$j_{\mu}$ to $j_{\sigma(\mu)}$. So the morphism $\coprod j_{\mu}$
descends to $\mathbb{F}_{p}$.

Suppose we have two isogenies \[
\rho:(E_{k}^{g},\mu_{k})\to(A,\lambda)\quad\textrm{and}\quad\rho':(E_{k}^{g},\mu_{k}')\to(A,\lambda)\]
with the same target, so that the corresponding points in $\coprod
X(\mathfrak{r}(w))$ map to the same point in $\mathcal{A}_{g}$. To
show that $\coprod j_{\mu}$ gives an isomorphism \[
\coprod_{\mu\in\Lambda_{g,c}}[X(\mathfrak{r}(w))/\mathrm{Aut}(E^{g},\mu)]\stackrel{\sim}{\to}\mathcal{S}_{w},\]
we must show that there is an automorphism $\gamma$ of $E^{g}$
such that $\gamma^{*}\mu=\mu'$.

By \cite{Ha07} proposition 3.1.5 both $\rho$ and $\rho'$ are
minimal isogenies in the terminology of \cite{LO98} section 1.8
(here we really use that $w$ is in $^{I}W_{g}^{(c)}$ and not just
$^{I}W_{g}^{[c]}$). Minimal isogenies are unique up to
isomorphism, so there is an isomorphism $\gamma$ of $E^{g}$ such
that $\rho'=\rho\circ\gamma$ and \[
\gamma^{*}\mu=\gamma^{*}(\rho^{*}\lambda)=(\rho')^{*}\lambda=\mu'\]
as required.

\section{\label{sec:Number-of-components}Number of components}

In this section we show using corollary
\ref{cor:properties-fine-DL-var} that $X(\mathfrak{r}(w))$ is
irreducible for $w$ in $^{I}W_{g}^{(c)}$. When we combine this
with theorem \ref{thm:Refinement-Harashita}, we get corollary
\ref{cor:Number-Of-Components}.

Let $W_{c}$ be the Weyl group of $\mathrm{Sp}_{2c}$ and let
$S_{c}=\{ s_{1},\dots,s_{c}\}\subset W_{c}$ be the set of
reflections in simple roots, where \[ s_{c}=(c,c+1)\textrm{ and
}s_{i}=(i,i+1)(2c-i,2c+1-i)\textrm{ for }i=1,\dots,c-1.\] For
$w\in W_{c}$ let $S_{c}(w)$ consist of those elements of $S_{c}$
that occur in a reduced expression for $w$. This set is
independent of the choice of a reduced expression. Also, $w$ is in
the subgroup $W_{c,J}$ generated by $J\subset S_{c}$ if and only
if $S_{c}(w)\subset J$ (see \cite{Bo02} chapter 4, section 8,
proposition 7 and corollary 1).

\begin{lem}
If $w$ is in $^{I}W_{g}^{(c)}$, then
$S_{c}(\mathfrak{r}(w))=S_{c}$.
\end{lem}
\begin{proof}
Example 3.6 in \cite{Mo01} tells us that any $w'\in{}^{I}W_{c}$
has a reduced expression of the form\[ w'=(s_{c}s_{c-1}\dots
s_{i_{l}})(s_{c}s_{c-1}\dots s_{i_{l-1}})\dots(s_{c}s_{c-1}\dots
s_{i_{1}})\] for some $i_{1}<\dots<i_{l}$. Hence $S_{c}(w')=\{
s_{i},s_{i+1},s_{i+2},\dots,s_{c}\}$ with $i=i_{1}$. So we are
done if we show that $s_{1}$ is in $S_{c}(\mathfrak{r}(w))$.

Note that $\mathfrak{r}(W_{g}^{[c-1]})=W_{c}^{[c-1]}$ is the
subgroup generated by $s_{2},s_{3},\dots,s_{c}$ and
$\mathfrak{r}(w)$ is not in this subgroup, since $w$ is in
$W_{g}^{(c)}=W_{g}^{[c]}-W_{g}^{[c-1]}$. So
$S_{c}(\mathfrak{r}(w))$ is not contained in $\{
s_{2},s_{3},\dots,s_{c}\}$ and it must contain $s_{1}$.
\end{proof}
Let $w$ be in $^{I}W_{g}^{(c)}$ and suppose that
$W_{I_{\infty}}\mathfrak{r}(w)$ is contained in $W_{c,J}$ for some
$J\subseteq S_{c}$ (with $I_{\infty}$ as in section
\ref{sub:Fine-Deligne-Lusztig-varieties}). Corollary
\ref{cor:properties-fine-DL-var} says that $X(\mathfrak{r}(w))$ is
irreducible if this implies that $J=S_{c}$. But $W_{c,J}$ contains
$\mathfrak{r}(w)$. So $J$ contains $S_{c}(\mathfrak{r}(w))$ and by
the lemma $S_{c}(\mathfrak{r}(w))=S_{c}$ .

\end{document}